\documentclass[10pt]{articlefedericosmallerlabels}

\usepackage{amsfonts}
\usepackage{amsmath}
\usepackage{amssymb}
\usepackage{amsthm}
\usepackage{upref}
\usepackage[colorlinks,final,backref=page,hyperindex]{hyperref}
\usepackage{mathtools}
\usepackage{float}

\usepackage{bbm}

\usepackage{xcolor}

\definecolor{darkgreen}{rgb}{0,.5,0}
\definecolor{brown}{rgb}{0.5,0.3,0}

\usepackage{url}
\usepackage{amsxtra}
\usepackage[final]{rotating}
\usepackage{array}
\usepackage[all]{xy}
\SelectTips{cm}{}
\usepackage{epic}

\usepackage{multicol}
\usepackage{caption}

\newtheorem{theorem}{Theorem}
\newtheorem{problem}[theorem]{Problem}

\theoremstyle{plain}

\newtheorem{definition}[theorem]{Definition}

\newtheorem*{axiom}{Axiom}

\newcommand{\C}{\mathcal{C}}
\renewcommand{\S}{\mathcal{S}}
\newcommand{\T}{\mathcal{T}}

\begin{document}
\title{CAT(0) geometry, robots, and society.}
\author{
\textsf{Federico Ardila--Mantilla\footnote{\noindent {Professor of Mathematics, San Francisco State University. Profesor Adjunto, Universidad de Los Andes. \texttt{federico@sfsu.edu}.
This work was supported by NSF CAREER grant DMS-0956178, grants DMS-0801075, DMS-1600609, DMS-1855610, Simons Fellowship 613384, and NIH grant 5UL1GM118985-03 awarded to SF BUILD.}}}}
\date{}

\maketitle

\begin{center}
\includegraphics[width=6in]{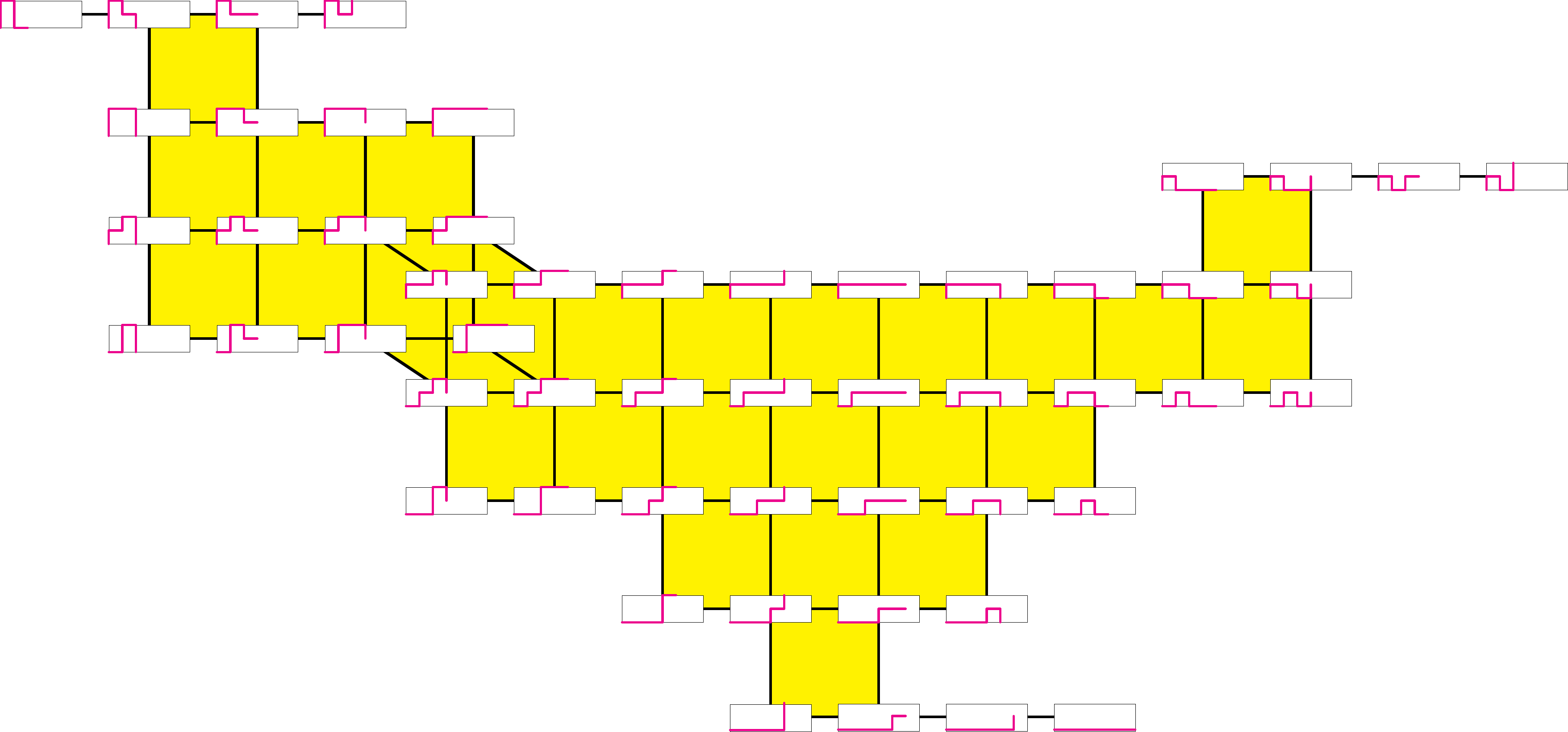}
\end{center}

\begin{multicols}{2}

\section{Moving objects optimally}

There are many situations in which an object can be moved using certain prescribed rules, and many reasons -- pure and applied -- to solve  the following problem.

\begin{problem}\label{prob:moverobot}
Move an object optimally from one given position to another.
\end{problem}	

\noindent Without a good idea, this is usually very hard to do.

When we are in a city we do not know well, trying to get from one location to another quickly, most of us will consult a map of the city to plan our route. This simple, powerful idea is the root of a very useful approach to Problem \ref{prob:moverobot}:
We build and understand the ``map of possibilities", which keeps track of all possible positions of the object; we call it the \emph{configuration space}. 
This idea is pervasive in many fields of mathematics, which call such maps \emph{moduli spaces}, \emph{parameter spaces}, or \emph{state complexes}.

This article seeks to explain that, for many objects that move discretely, the resulting ``map of possibilities" is a \emph{CAT(0) cubical complex}: a space of non-positive curvature made of unit cubes. When this is the case, we can use ideas from geometric group theory and combinatorics to solve Problem \ref{prob:moverobot}.
 
 \columnbreak
 
This approach is applicable to many different settings; but to keep the discussion concrete, we focus on the following specific example. For precise statements, see Section \ref{sec:movingrobots} and Theorems  \ref{thm:algorithm} and \ref{th:CAT(0)2}.

\begin{center}
	\includegraphics[height = 1.5cm]{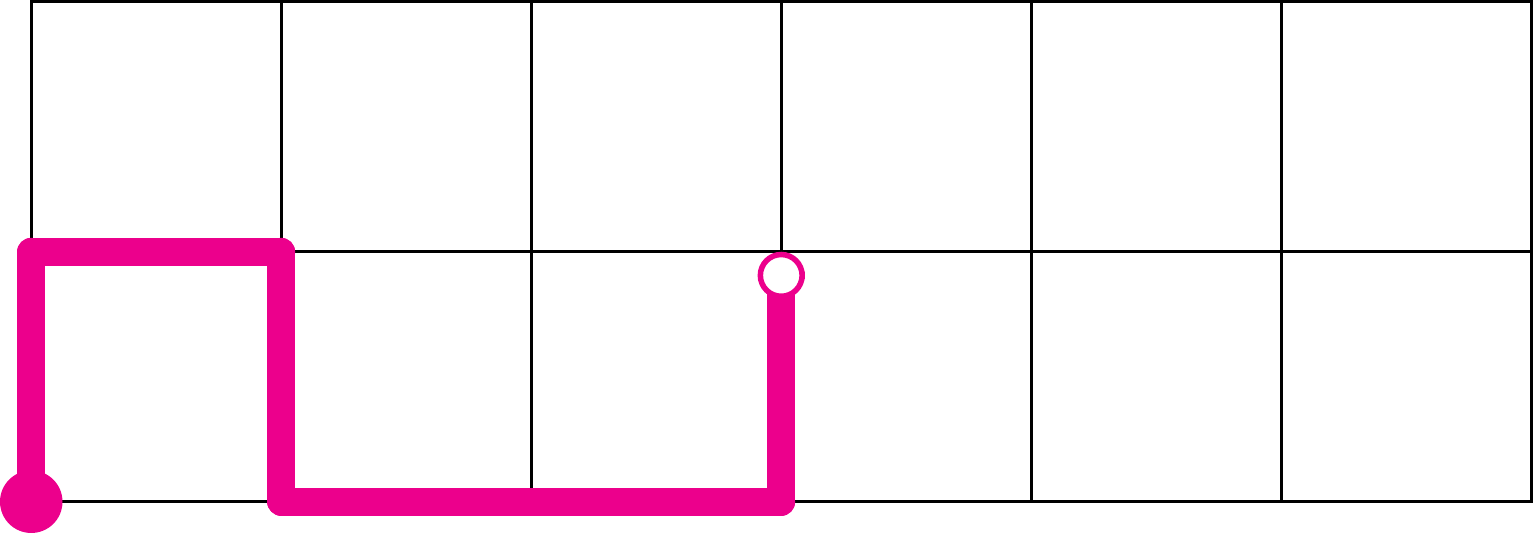} \qquad
	\captionof{figure}{\label{fig:armingrid} A pinned down robotic arm of length 6 in a tunnel of height 2. The figure above shows its configuration space.}
\end{center}

\begin{theorem} \cite{ArdilaBakerYatchak, ArdilaBastidasCeballosGuo}\label{th:CAT(0)}
The configuration space of a 2-D pinned down robotic arm in a rectangular tunnel is a $\mathrm{CAT}(0)$ cubical complex. Therefore there is an explicit algorithm to move this robotic arm optimally from one given position to another.
\end{theorem}


\section{Black Lives Matter}

On July 4, 2016 we finished the implementation of our algorithm to move a discrete robotic arm. Three days later, seemingly for the first time in history, US police used a robot to kill an American citizen. Now, whenever I present this research, I also discuss this action. 

The day started with several peaceful Black Lives Matter protests across the US, condemning the violence disproportionately inflicted on Black communities by the American state. These particular protests were prompted by the shootings of Alton Sterling and Philando Castile by police officers in Minnesota and Louisiana.

In Dallas, TX, as the protest was coming to an end, a sniper opened fire on the crowd, killing five police officers. Dallas Police initially misidentified a Black man -- the brother of one of the protest organizers -- as a suspect. They posted a photo of him on the internet and asked for help finding him. Fearing for his life, he turned himself in, and was quickly found innocent. 

A few hours later, police identified US Army veteran Micah Johnson as the main suspect. After a chase, a standoff, and failed negotiations, they used a robot to kill him, without due process of law.

The organizers of the protest condemned the sniper's actions, and police officials believe he acted alone. The robot that killed Johnson cost about \$150,000; police said that the arm of the robot was damaged, but still functional after the blast \cite{CNN}. The innocent man who was misidentified by the police continued to receive death threats for months afterwards.

Different people will have different opinions about the actions of the Dallas Police in this tragic event. 
What is certainly unhealthy is that the large majority of people I have spoken to have never heard of this incident. 

Our mathematical model of a robotic arm is very simplified, and probably far from direct applications, but the techniques developed here have the potential to make robotic operations cheaper and more efficient. We tell ourselves that mathematics and robotics are neutral tools, but our research is not independent from how it is applied. We arrive to mathematics and science searching for beauty, understanding, or applicability. When we discover the power that they carry, how do we proceed?

\begin{axiom} \cite{TodosCuentan} Mathematics is a powerful, malleable tool that can be shaped and used differently by various communities to serve their needs.
\end{axiom}

Who currently holds that power? How do we use it? Who funds it and for what ends? 
With whom do we share that power? Which communities benefit from it? Which are disproportionately harmed by it?

For me these are the hardest questions about this work, and the most important. 
The second goal in writing this article -- a central one for me -- is to invite myself, and its readers, to continue to look for answers that make sense to us.

\section{Moving robots}\label{sec:movingrobots}

We consider a discrete 2-D robotic arm $R_{m,n}$ of length $n$ moving in a rectangular tunnel of height $m$. The robot consists of $n$ links of unit length, attached sequentially, facing up, down, or right. Its base is affixed to the lower left corner of the tunnel, as shown in Figure~\ref{fig:armingrid} for $R_{2,6}$.

The robotic arm may move freely, as long as it doesn't collide with itself, using two kinds of local moves:

\noindent $\bullet$ \emph{Flipping a corner:} Two consecutive links facing different directions interchange directions.

\noindent $\bullet$ \emph{Rotating the end:} The last link of the robot rotates $90^\circ$ without intersecting itself.

\noindent
This is an example of a \emph{metamorphic robot} \cite{AbramsGhrist}.

\begin{center}
    \includegraphics[height = 1.5cm]{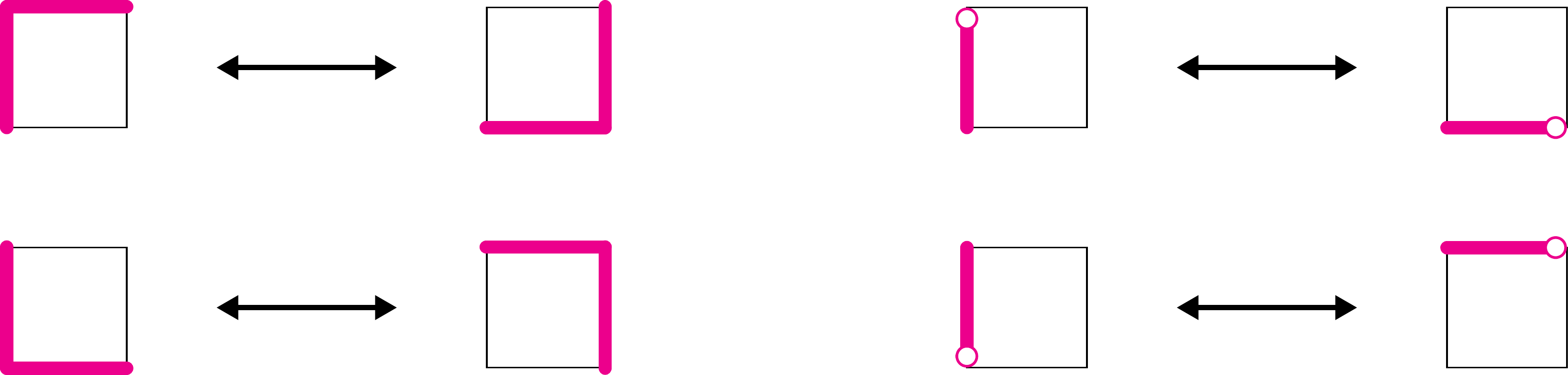}
	\captionof{figure}{The local moves of the robotic arm.}
	\label{fig:mov}
\end{center}	

\noindent 
How can we get the robot to navigate this tunnel efficiently?

Figure \ref{fig:twopositions} shows two positions of the robot; suppose we want to move it from one position to the other. By trial and error, one will not have too much difficulty in doing it. It is not at all clear, however, how one might do this in the most efficient way possible.

\begin{center}
	\includegraphics[height = 1cm]{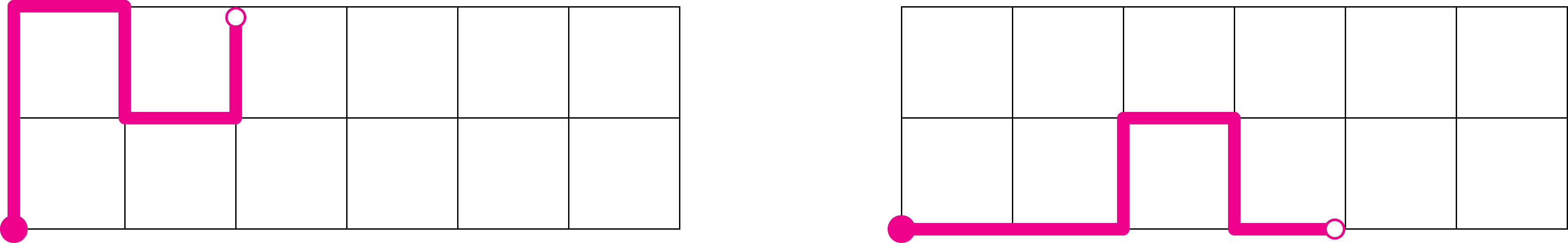}
	\captionof{figure}{\label{fig:twopositions} Two positions of the robot $R_{2,6}$.}
\end{center}

\section{Maps}

To answer this question, let us build the ``map of possibilities" of the robot. We begin with a  \emph{configuration graph}, which has a node for each position of the arm, and an edge for each local move between two positions. A small piece of this graph is shown below.

\begin{center}
\includegraphics[width=3in]{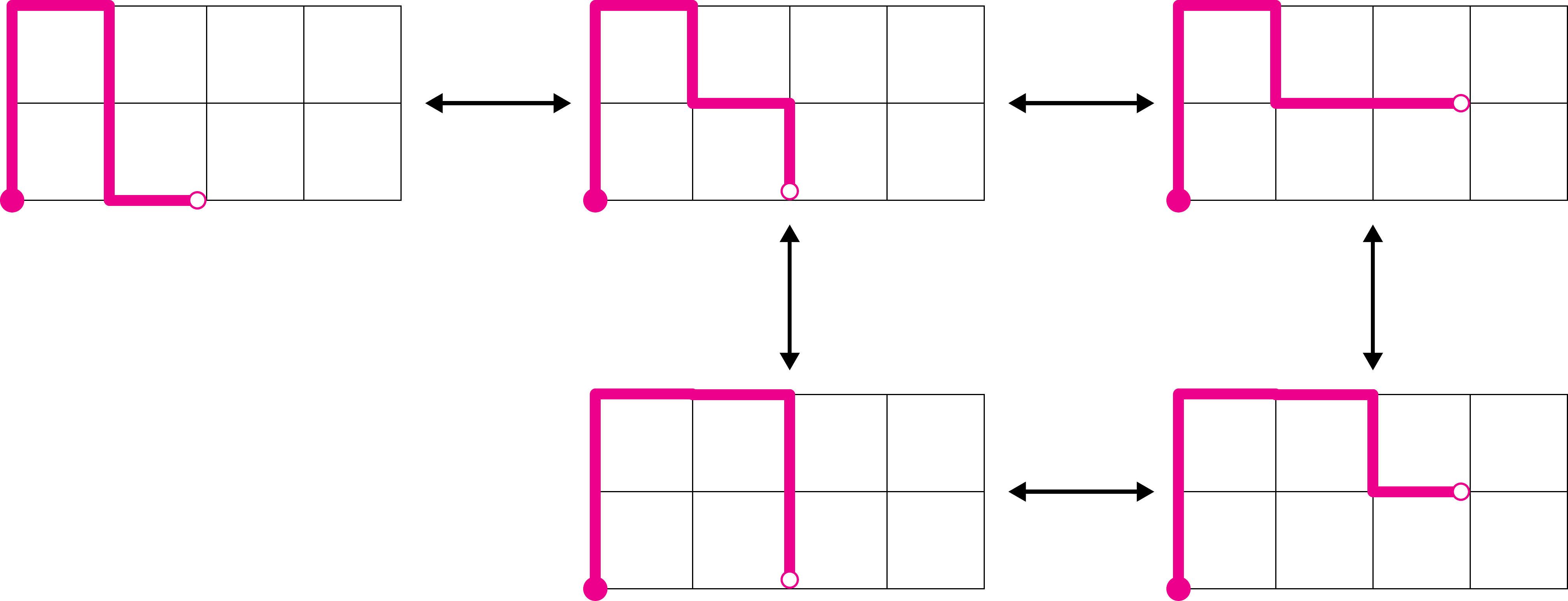}
	\captionof{figure}{\label{fig:thegraph} A part of the graph of possibilities of $R_{2,6}$.}
\end{center}

As we see in the figure of page 1, the resulting graph looks a bit like the map of downtown San Francisco or Bogot\'a, with many square blocks lined up neatly. Such a cycle of length $4$ arises whenever the robot is in a given position, and there are two moves A and B that do not interfere with each other: if we perform move A and then move B, the result is the same as if we perform move B and then move A; see for example the $4$-cycle of Figure \ref{fig:thegraph}. 
More generally, if the robot has $k$ moves that can be performed independently of each other, these moves result in (the skeleton of) a $k$-dimensional cube in the graph. 

This brings up an important point: If we wish to move the robot efficiently, we should let it perform various moves simultaneously. In the map, this corresponds to walking across the diagonal of the corresponding cube. Thus we construct the \emph{configuration space} of the robot, by filling in the $k$-cube corresponding to any $k$ moves that can be performed simultaneously, as illustrated in Figure \ref{fig:thecomplex}; compare with Figure \ref{fig:thegraph}.
The result is a \emph{cubical complex}, a space made of cubes that are glued face-to-face.

\begin{center}
\includegraphics[width=2in]{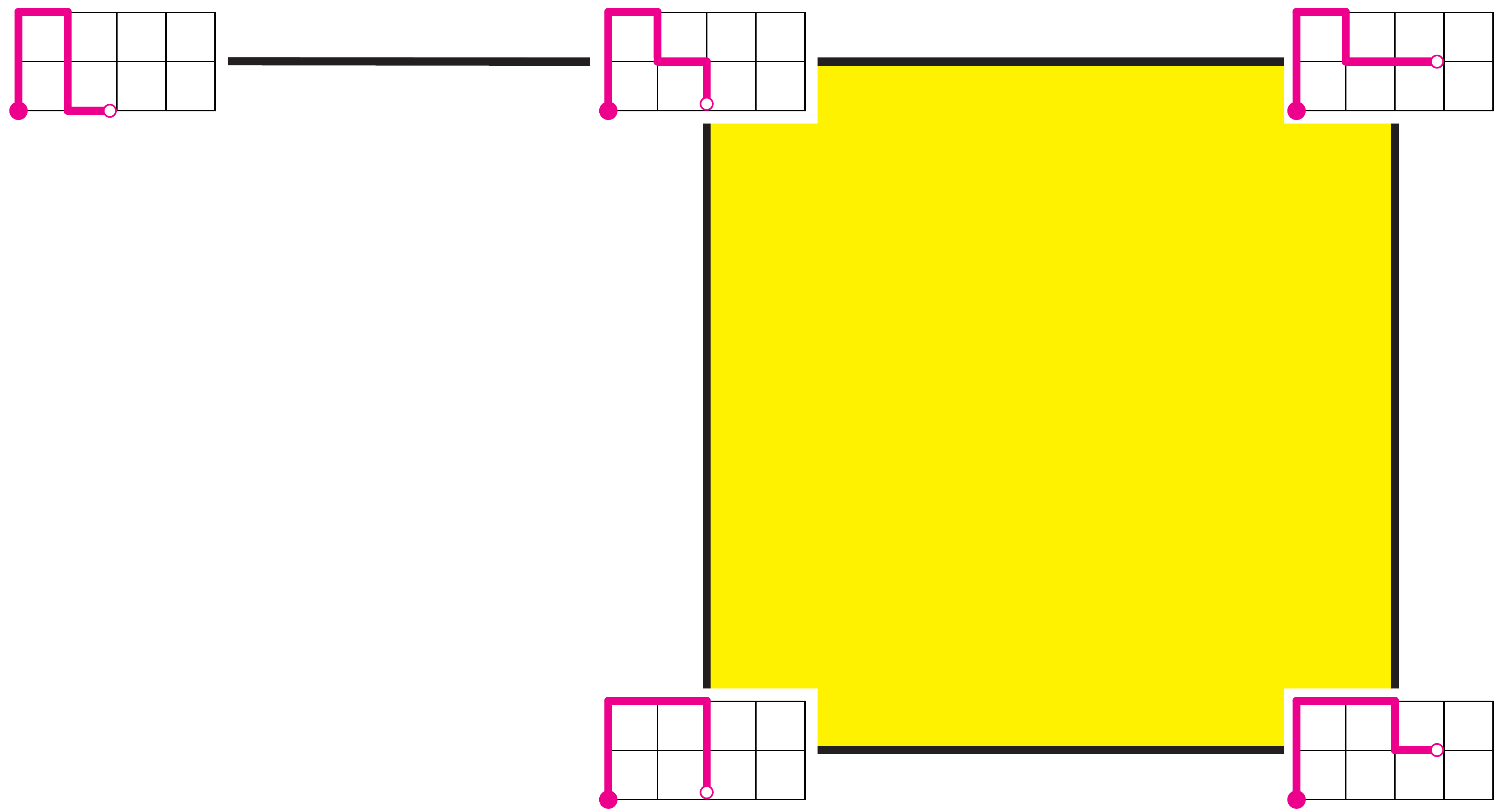}
	\captionof{figure}{\label{fig:thecomplex} A part of the configuration space of $R_{2,6}$.}
\end{center}

\begin{definition} \label{def:configspace} The \emph{configuration space} $\C(R)$ of the robotic arm $R$ is the cubical complex with:

$\bullet$ a vertex for each position of the robot,

$\bullet$ an edge for each local move between two positions, 

$\bullet$ a $k$-dimensional cube for each $k$-tuple of local moves that may be performed simultaneously.

\end{definition}

This definition applies much more generally to discrete situations that change according to local moves; see Section \ref{sec:examples} and \cite{AbramsGhrist, GhristPeterson}.

In our specific example, Figure \ref{fig:configspace} shows the configuration space of the robot $R_{2,6}$ of length $6$ in a tunnel of height $2$. It is now clear how to move between two positions efficiently: just follow the shortest path between them in the map!

\begin{center}
\includegraphics[width=3.5in]{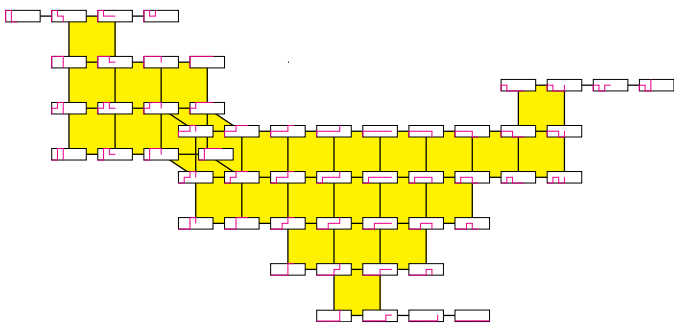}
\captionof{figure}{\label{fig:configspace}The configuration space of the robotic arm $R_{2,6}$.}
\end{center}

\section{What are we optimizing?}

Is it so clear, just looking at a map, what the optimal path will be? It depends on what we are trying to optimize. In San Francisco, with its beautifully steep hills, the best route between two points can be very different depending on whether one is driving, biking, walking, or taking public transportation. 
The same is true for the motion of a robot.

For the configuration spaces we are studying, there are at least three reasonable metrics: $\ell_1, \ell_2$, and $\ell_\infty$. In these metrics, the distance between points $x$ and $y$ in the same $d$-cube, say $[0,1]^d$, is
\[
\sqrt{\sum_{1 \leq i \leq d} (x_i-y_i)^2}, \qquad
\sum_{1 \leq i \leq d} |x_i-y_i|, \qquad
\max_{1 \leq i \leq d} |x_i-y_i|,
\]
respectively. 
Figure \ref{fig:shortest} shows the two positions of the robot of Figure \ref{fig:twopositions} in the configuration space, and shortest paths or \emph{geodesics} between them according to these metrics. 

\begin{center}
	\includegraphics[height = 4cm]{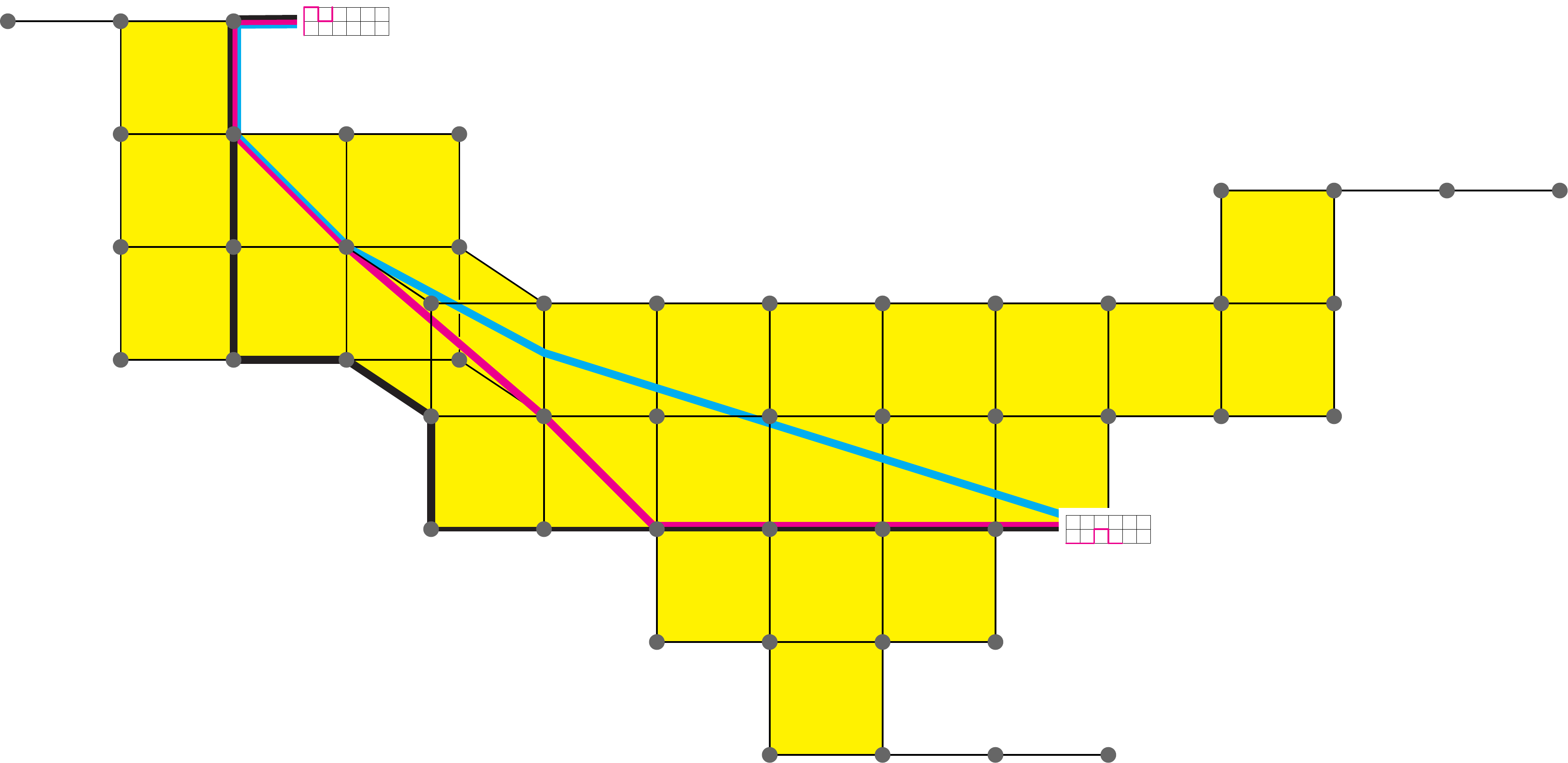}
	\captionof{figure}{\label{fig:shortest} Some paths between two points in the configuration space $\C(R_{2,6})$. 
	The black path is geodesic in the $\ell_1$ metric, the magenta path is geodesic in  $\ell_1$ and $\ell_\infty$, and the cyan path is geodesic in $\ell_1, \ell_2,$ and $\ell_\infty$.
	}

\end{center}


If each individual move has a ``cost" of $1$, then performing $d$ simultaneous moves -- which corresponds to crossing a $d$-cube -- costs $\sqrt{d}$, $d$, and $1$ in the metrics $\ell_2$, $\ell_1$, and $\ell_\infty$. Although the Euclidean metric is the most familiar, it seems unrealistic in this application; why should  two simultaneous moves cost $\sqrt 2$? For the applications we have in mind, the $\ell_1$ and $\ell_\infty$ metrics are reasonable models for the cost and the time of motion:

\smallskip
\noindent
\emph{Cost ($\ell_1$)}: We perform one move at a time; there is no cost benefit to making moves simultaneously.

\smallskip
\noindent
\emph{Time ($\ell_\infty$)}: We may perform several moves at a time, causing no extra delay.

\smallskip
These two metrics, studied in \cite{ArdilaBakerYatchak, ArdilaBastidasCeballosGuo}, will be the ones that concern us in this paper. The Euclidean metric, which is useful in other contexts and significantly harder to analyze, is studied in \cite{ArdilaOwenSullivant, Hayashi}.

\section{Morning routine}

I write this while on sabbatical in a foreign city. 
Being the coffee enthusiast that I am, I carefully study a map several mornings in a row, struggling to find the best cafe on my way from home to my office. One morning, amused, my partner May-Li stops me on the way out:

\medskip
\noindent
-- Fede, you know you don't \textbf{always} have to take a geodesic, right?

\medskip

Perhaps, instead of the most efficient paths, we should be looking for the most pleasant, or the greenest, or the most surprising, or the most beautiful.

\section{CAT(0) cubical complexes}

Our two most relevant algorithmic results are the explicit construction of cheapest ($\ell_{1}$) and fastest ($\ell_\infty$) paths in the configuration space of the robot arm $R_{m,n}$. Still, the Euclidean metric ($\ell_2$) turns out to play a very important role as well. Most configuration spaces that interest us exhibit non-positive curvature with respect to the Euclidean metric, and this fact is central in our construction of shortest paths in the cost and time metrics.

Let us consider a \emph{geodesic metric space} $(X,d)$, where any two points $x$ and $y$ can be joined by a unique shortest path of length $d(x,y)$; such a path is known as a \emph{geodesic}. 
Let $T$ be a triangle in $X$ whose sides are geodesics of lengths $a,b,c$, and let $T'$ be the triangle with the same sidelengths in the plane. For any geodesic chord of length $d$ connecting two points on the boundary of $T$, there is a comparison chord between the corresponding two points on the boundary of $T'$, say of length $d'$. If $d \leq d'$ for any such chord in $T$, we say that triangle $T$ is \emph{at least as thin as a Euclidean triangle}.

\begin{center}
\includegraphics[scale=.7]{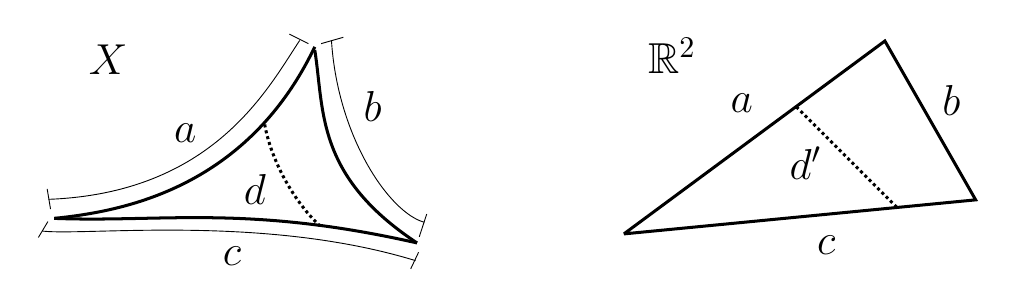}
\captionof{figure}{A chord in a triangle in $X$, and the corresponding chord in the comparison triangle in $\mathbb{R}^2$. The triangle in $X$ is \emph{at least as thin as a Euclidean triangle} if $d \leq d'$ for all such chords. \label{fig:thintriangle}}
\end{center}

\begin{definition}\label{def:CAT(0)}
A metric space $X$ is $\mathrm{CAT}(0)$  if:

\noindent $\bullet$ 
between any two points there is a unique geodesic, and 

\noindent $\bullet$
every triangle is at least as thin as a Euclidean triangle.
 \end{definition}

A (finite) \emph{cubical complex} is a connected space obtained by gluing finitely many cubes of various dimensions along their faces. 
We regard it as a metric space with the Euclidean metric on each cube; all cubes necessarily have the same side length. Cubical complexes are flat inside each cube, but they can have curvature where cubes are glued together, for example, by attaching three or five squares around a common vertex (obtaining positive and negative curvature, respectively), as shown in Figure \ref{fig:fatthin}. 
We invite the reader to check that the triangles in the left and right panel of Figure \ref{fig:fatthin} are thinner and not thinner than a Euclidean triangle, respectively.

\begin{center}
\includegraphics[height=1.2in]{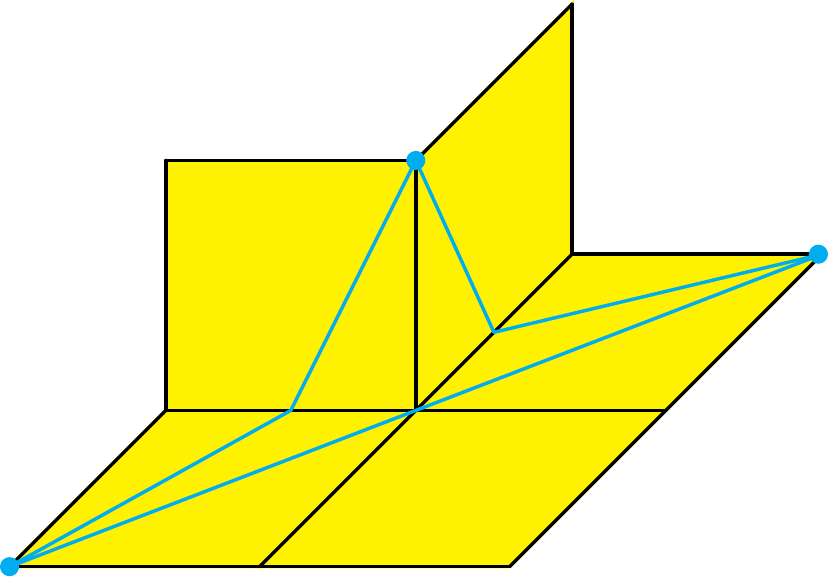} \qquad
\includegraphics[height=.9in]{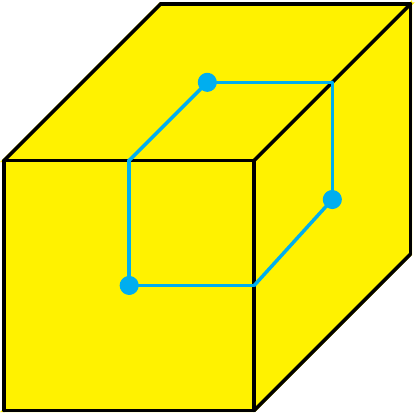} 
\captionof{figure}{A CAT(0) and a non-CAT(0) cubical complex.}
\label{fig:fatthin}
\end{center}

We have the following general theorem.

\begin{theorem} \label{thm:algorithm} \cite{AbramsGhrist, ArdilaBakerYatchak, ArdilaOwenSullivant, Hayashi, NibloReeves} 
Given two points $x$ and $y$ in a CAT(0) cubical complex, there are algorithms to find a geodesic from $x$ to $y$ in the Euclidean ($\ell_2$), cost ($\ell_1$), and time ($\ell_\infty$) metrics.
\end{theorem}

Thus a robot with a CAT(0) configuration space is easier to control: 
we have a procedure that automatically moves it optimally. 
We will see that this is the case for the robotic arm $R_{m,n}$.

\section{How do we proceed?}

Once I began to feel that this work, which started out in ``pure" mathematics, could actually have real-life applications, I started getting anxious and selective about who I discussed it with. It is a strange feeling, to discover something you really like, and yet to hope that not too many people find out about it. 
When I was invited to write this article, I felt conflicted. 
I knew I did not want to only discuss the mathematics, but I am much less comfortable writing outside of the shared imaginary world of mathematicians, where we believe we know right from wrong. Still, I know it is important to listen, learn, discuss, and even write from this place of discomfort. 

How should I tell this story? Should I do it at all? I have turned to many friends, colleagues, and students for their wisdom and advice.

\medskip

Mario Sanchez, who thinks deeply and critically about the culture of mathematics and philosophy in our society, is wary of mathematical fashions: What if it becomes trendy for mathematicians to start working on optimizing robots, but not to think about what is being optimized, or whom that optimization benefits? 
He tells me, with his quiet intensity: ``If you're worried that your paper might have this effect, you should probably emphasize the human question pretty strongly."

Laura Escobar just returned from a yoga retreat in Champaign-Urbana where a scholar of Indian literature taught them the story of Arjuna, a young warrior about to enter a rightful battle against members of his own family. Deeply conflicted about the great violence that will ensue, he turns to Krishna for advice.
Oversimplifying his reply, Krishna says: ``One should not abandon duties born of one?s nature, even if one sees defects in them. It is your duty as a warrior to uphold the Dharma, take action, and fight."
With her usual thoughtful laugh, she tells me about the distressed reactions of her peace-loving yoga classmates. Laura and I grew up in the middle of Colombia's 60--year old civil war, which has killed more than 215,000 civilians and 45,000 combatants and has displaced more than 15\% of the country's population \cite{CNMH, UN}; it is hard for us to understand Krishna's advice as well.\footnote{We later learn that Robert Oppenheimer quoted Krishna when he and his team detonated the first nuclear bomb.} So we go to the bookstore and buy matching copies of the Bhagavad Gita.

Many of my friends who do not work in science are surprised by the lack of structural and institutional resources. They ask me: If a mathematician or a scientist is trying to understand or have some control over the societal impact of their knowledge and their expertise, what organizations they can turn to for support? I have been asking this question to many people. I have not found one, but I am collecting resources. 
Interdisciplinary organizations like the Union of Concerned Scientists, Science for the People, Data for Black Lives, and sections of the American Association for the Advancement of Science seek to use science to improve people's lives and advance social justice. 
Our colleagues in departments of Science, Technology, and Society, Public Policy, History, Philosophy, and Ethnic Studies have been studying these issues for decades, even centuries. 
This has often taken place too far from science departments, and it must be said that my generation of scientists largely looked down on these disciplines as unrigorous, uninteresting, or unimportant. 
Governments, companies, and professional organizations assemble Ethics Committees, usually separate from their main operations, and give them little to no decision-making power. 

How do we make these considerations an integral part of the practice and application of science? I am encouraged to see that the new generation of scientists understands their urgent role in society much more clearly than we do. 

May-Li Khoe, whom I can always trust to be wise and direct, asks me: If you tell me that this model of mapping possibilities could be applicable in many areas, and you don't trust the organizations that build the most powerful robots, why don't you find other applications? She's right. I'm looking.

\section{Examples. \label{sec:examples}}

Just like any other cultural practice, mathematics respects none of the artificial boundaries that we sometimes draw, in an attempt to understand it and control it. This is evident for CAT(0) cubical complexes, a family of objects which appears in many seemingly disparate parts of (mathematical) nature. Let us discuss three sources of examples; each one raises different kinds of questions and offers valuable tools that have directly shaped this investigation.

\smallskip

\noindent \textsf{\textbf{Geometric group theory.}}
This project was born in geometric group theory, which studies groups by analyzing how they act on geometric spaces. Gromov's pioneering work in this field  \cite{Gromov} led to the systematic study of CAT(0) cubical complexes. A concrete source of examples is due to Davis \cite{Davis}.

A \emph{right-angled Coxeter group} $X(G)$ is given by generators of order 2 and some commuting relations between them; we encode the generators and commuting pairs in a graph $G$. For example, the graph of Figure \ref{fig:davis}.a encodes the group generated by $a,b,c,d$ with relations $a^2=b^2=c^2=d^2 = 1,$
 $ ab=ba, ac=ca, bc=cb, cd=dc $.

The \emph{Cayley graph} has a vertex for each element of $X(G)$ and an edge between $g$ and $gs$ for each group element $g$ and generator $s$. This graph is the skeleton of a CAT(0) cube complex that $G$ acts on, called the \emph{Davis complex} $\S(G)$. It is illustrated in Figure \ref{fig:davis}.b. One can then use the geometry of $\S(G)$ to derive algebraic properties of $X(G)$. For example, one can easily solve the \emph{word problem} for this group: given a word in the generators, determine whether it equals the identity. This problem is undecidable for general groups.

\begin{center}
\includegraphics[height=.92in]{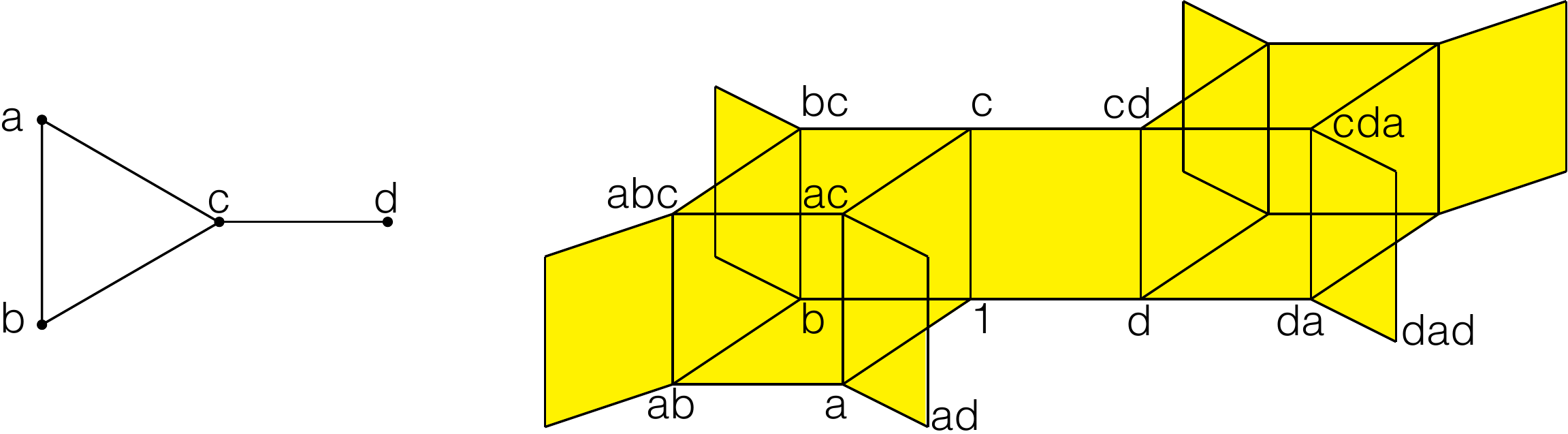} 
\captionof{figure}{a. A graph $G$ determining a right-angled Coxeter group $X(G)$, and b. part of its Davis complex $\S(G)$.\label{fig:davis}}
\end{center}

\smallskip

\noindent \textsf{\textbf{Phylogenetic trees.}} A central problem in phylogenetics is the following: given $n$ species, determine the most likely evolutionary tree that led to them. There are many ways of measuring how different two species are\footnote{We should approach them thoughtfully and critically; see Section \ref{sec:phylo}.}; but if we are given the ${n \choose 2}$ pairwise distances between the species, how do we construct the tree that most closely fits that data?

\begin{center}
\includegraphics[height=3.5in]{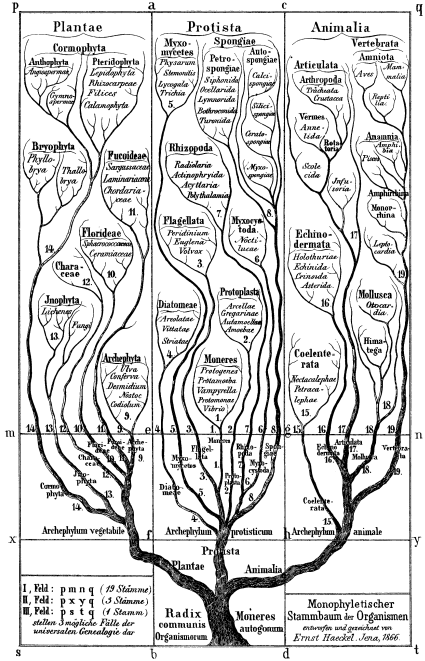} 
\captionof{figure}{\label{fig:tree}Ernst Haeckel's tree of life (1866).}
\end{center}

Billera, Holmes, and Vogtmann \cite{BilleraHolmesVogtmann} approached this problem by constructing the space of all possibilities: the \emph{space of trees} $\T_n$. 
Remarkably, they proved that the space of trees $\T_n$ is a CAT(0) cube complex. In particular, since it has unique geodesics, we can measure the distance between two trees, or find the average tree between them. This can be very helpful in applications: if 10 different algorithms propose 10 different phylogenetic trees, we can detect which proposed trees are close to each other, detect outlier proposals that seem unlikely, or find the average between different proposals. Owen and Provan showed how to do this in polynomial time \cite{OwenProvan}.

\begin{center}
\includegraphics[height=1.4in]{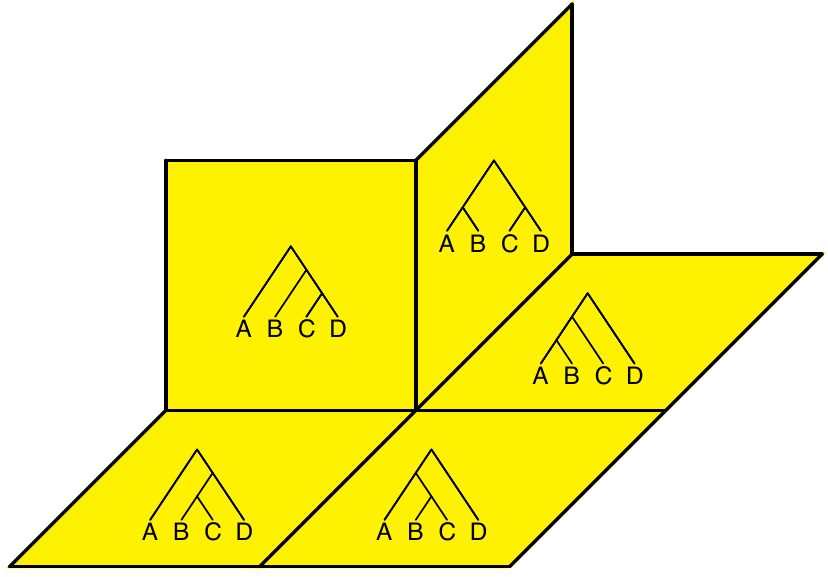} 
\captionof{figure}{Five of the 15 squares in the space of trees $\T_4$.}
\end{center}

These results made us wonder whether one can similarly construct $\ell_2$-optimal paths in any CAT(0) cube complex. New complications arise, but it \textbf{is} possible. \cite{ArdilaOwenSullivant, Hayashi}.

\smallskip

\noindent \textsf{\textbf{Discrete systems: reconfiguration.}}
Abrams, Ghrist, and Peterson introduced \emph{reconfigurable systems} in  \cite{AbramsGhrist, GhristPeterson}. This very general framework models discrete objects that change according to local moves, keeping track of which pairs of moves can be carried out simultaneously. Examples include discrete metamorphic robots moving around a space, particles moving around a graph without colliding, domino tilings of a region changing by \emph{flips} \includegraphics[height=.15in]{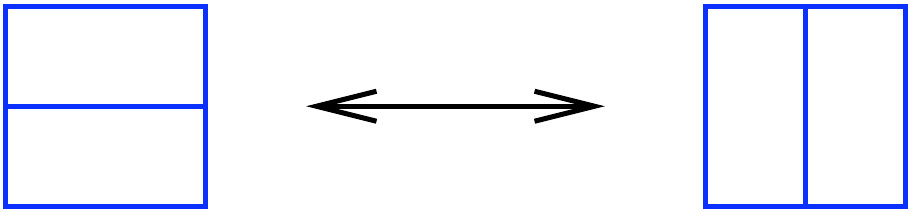}, and reduced words in the symmetric group changing by commutation moves $s_is_j \leftrightarrow s_js_i$ for $|i - j| \geq 2$ and braid moves $s_is_{i+1}s_i \leftrightarrow s_{i+1}s_is_{i+1}$.

Definition \ref{def:configspace} associates a configuration space to any reconfigurable system. 
Such a configuration space is always \textbf{locally} CAT(0). It is often \textbf{globally} CAT(0), and when that happens Theorem \ref{thm:algorithm} applies, allowing us to move our objects optimally.

\section{Why do we map?}\label{sec:phylo}

Math historian Michael Barany points out to me that, struck by the aesthetic beauty of the tree of life shown in Figure \ref{fig:tree}, I failed to notice another map that Haeckel drew: a hierarchical tree of nine human groups -- which he regarded as different species -- showing their supposed evolutionary distance from the ape-man. 
Modern biology shows this has no scientific validity, and furthermore, that there is no genetic basis for the concept of race. Haeckel's work is just one sample of the deep historical ties between phylogenetics and scientific racism, and between mapmaking and domination.

\begin{quote}
If we map from a different -- an \emph{other} -- point of view [...]
then mapping becomes a process of getting to know, connect, bring closer together in relation, remember, and interpret. 

-- Sandra Alvarez \cite{Alvarez}
\end{quote}

\section{Characterizations}

How does one determine whether a given space is CAT(0)? We surely do not want to follow Definition \ref{def:CAT(0)} and check whether \textbf{every} triangle is at least as thin as a Euclidean triangle; this is not easy to do, even for an example as small as Figure \ref{fig:fatthin}. Fortunately, this becomes much easier when the space in question is a cubical complex. In this case, Gromov showed that the CAT(0) property -- a subtle metric condition -- can be rephrased entirely in terms of topology and combinatorics; no measuring is necessary! 

To state this, we recall two definitions. A space $X$ is \emph{simply connected} if there is a path between any two points, and every loop can be contracted to a point. If $v$ is a vertex of a cubical complex $X$, then the \emph{link} of $v$ in $X$ is the simplicial complex one obtains by intersecting $X$ with a small sphere centered at $v$. A simplicial complex $\Delta$ is \emph{flag} if it has no empty simplices: if $A$ is a set of vertices and every pair of vertices in $A$ is connected by an edge in $\Delta$, then $A$ is a simplex in $\Delta$.

\begin{theorem} \cite{Gromov}
A cubical complex $X$ is CAT(0) if and only if:

$\bullet$ $X$ is simply connected, and

$\bullet$ the link of every vertex in $X$ is flag.
\end{theorem} 

In fact, one can also do without the topology: there is an entirely combinatorial characterization of CAT(0) cubical complexes. This is originally due to Sageev and Roller, and we rediscovered it in \cite{ArdilaOwenSullivant} in a different formulation that is more convenient for our purposes.
Let a \emph{pointed cubical complex} be a cubical complex with a distinguished vertex.

\begin{definition} \cite{ArdilaOwenSullivant, Winskel}
A \emph{poset with inconsistent pairs (PIP)} $(P, \leq, \nleftrightarrow)$  is a poset $(P, \leq)$ together with a collection of \emph{inconsistent pairs}, denoted $p \nleftrightarrow q$ for $p \neq q \in P$, that is closed under $\leq$; that is, 

if $p \nleftrightarrow q$ and $p \leq p'$, $q \leq q'$, then 
$p' \nleftrightarrow q'$.
\end{definition}

\noindent PIPs are also known as \emph{prime event structures} in the computer science literature \cite{Winskel}. The \emph{Hasse diagram} of a PIP $(P, \leq, \nleftrightarrow)$ shows graphically the minimal relations that define it. It has a dot for each element of $P$, a solid line from $p$ upward to $q$ whenever $p<q$ and there is no $r$ with $p<q<r$, and a dotted line between $p$ and $q$ whenever $p \nleftrightarrow q$ and there are no $r \leq p$ and $s \leq q$ such that $r \nleftrightarrow s$.

\begin{center}
\includegraphics[height=1in]{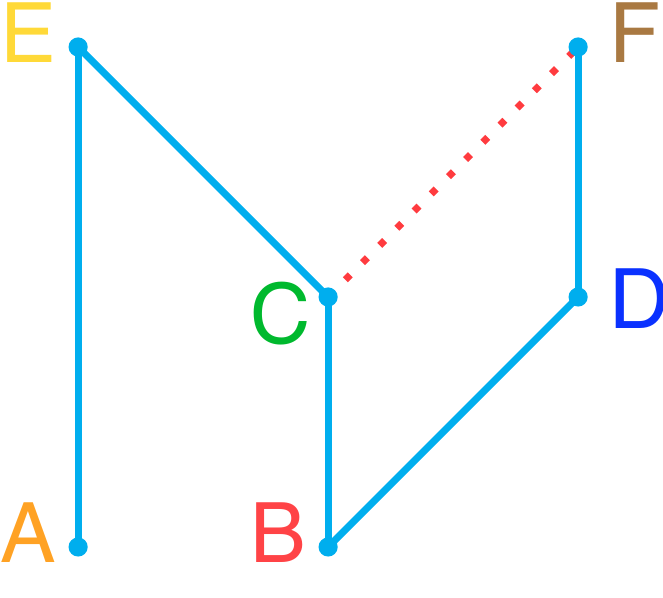}  \qquad \qquad 
\captionof{figure}{The Hasse diagram of a PIP: solid lines represent the poset, and dotted lines represent the (minimal) inconsistent pairs. Notice that $C \nleftrightarrow F$ implies $E \nleftrightarrow F$.\label{fig:PIP}}
\end{center}

\begin{theorem} \label{thm:PIP}
\cite{ArdilaOwenSullivant, Roller, Sageev} Pointed CAT(0) cube complexes are in bijection with posets with inconsistent pairs (PIPs).
\end{theorem}

This rediscovery was motivated by the observation that CAT(0) cubical complexes look very much like distributive lattices. In fact, Theorem \ref{thm:PIP} is an analog of Birkhoff's representation theorem, which gives a bijection between distributive lattices and posets. The proof is subtle and relies heavily on Sageev's work \cite{Sageev}, but the bijection is easy and useful to describe:

\smallskip

\noindent
\textsf{Pointed CAT(0) cubical complex $\mapsto$ PIP}:
Let $(X,v)$ be a CAT(0) cubical complex $X$ rooted at vertex $v$. Every $d$-cube in $X$ has $d$ hyperplanes that bisect its edges. Whenever two cubes share an edge, let us glue the two hyperplanes bisecting it. The result is a system of \emph{hyperplanes} associated to $X$ \cite{Sageev}. Figure \ref{fig:walls} shows an example. 

The PIP corresponding to $(X,v)$ keeps track of how one can navigate $X$ starting from $v$. 
The elements of the corresponding PIP are the hyperplanes. We declare $H < I$ if, starting from $v$, one must cross $H$ before crossing $I$. We declare $H \nleftrightarrow I$ if, starting from $v$, one cannot cross both $H$ and $I$ without backtracking. Remarkably, the simple combinatorial information stored in this PIP is enough to recover the space $(X,v)$.

\begin{center}
\includegraphics[height=1.5in]{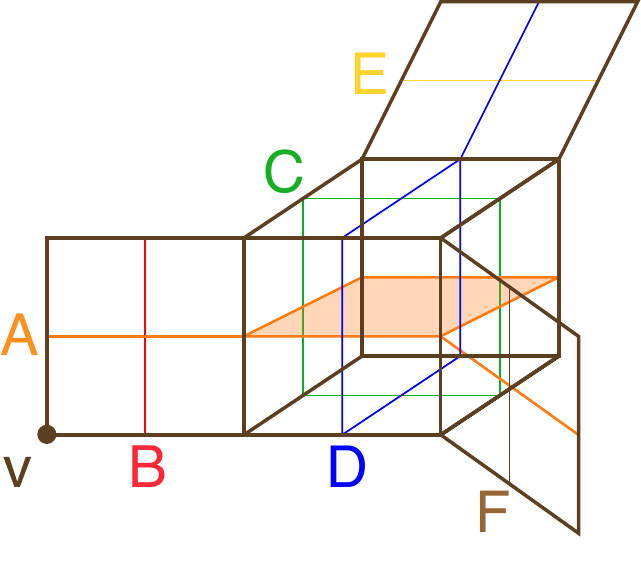}  \qquad \qquad 
\captionof{figure}{A rooted CAT(0) cubical complex with six hyperplanes. Its PIP is shown in Figure \ref{fig:PIP}.\label{fig:walls}}
\end{center}

\smallskip

\noindent 
\textsf{PIP $\longmapsto$ rooted CAT(0) cubical complex}: Let $P$ be a PIP.
An \emph{order ideal} of $P$ is a subset $I$ closed under $<$; that is, if $x<y$ and $y \in I$ then $x \in I$. We say that $I$ is \emph{consistent} if it contains no inconsistent pair.

The vertices of the corresponding CAT(0) cubical complex $X(P)$ correspond to the consistent order ideals of $P$. Two vertices are connected if their ideals differ by a single element. Then we fill in all cubes whose edges are in this graph. The root is the vertex corresponding to the empty order ideal.

We invite the reader to verify that the PIP of Figure \ref{fig:PIP} corresponds to the rooted complex of Figure \ref{fig:walls}.

\medskip

Theorem \ref{thm:PIP} provides a completely combinatorial way of proving that a cubical complex is CAT(0): one simply needs to identify the corresponding PIP!

\section{Remote controls and geodesics}

Intuitively, we think of the PIP $P$ as a ``remote control" to help an imaginary particle navigate the corresponding CAT(0) cubical complex $X$. If the particle is at a vertex of $X$, there is a corresponding consistent order ideal $I$ of $P$. The hyperplanes that the particle can cross are the maximal elements of $I$ and the minimal elements of $P-I$ consistent with $I$. We can then press the $i$th ``button" of $P$ if we want the point to cross hyperplane $i$.

This point of view is powerful because in practical applications, the configuration space $X$ is usually very large, high dimensional, and combinatorially complicated, whereas the remote control $P$ is much smaller and can be constructed in some cases of interest.

Theorem \ref{thm:algorithm} provides algorithms to move optimally between any two points in a CAT(0) cubical complex in the $\ell_1, \ell_2,$ and $\ell_\infty$ metrics. We sketch the proof in the cases that are relevant here: in the $\ell_1$ and $\ell_\infty$ metrics, where the two points $v$ and $w$ are vertices.

\begin{proof}[Sketch of Proof of Theorem \ref{thm:algorithm}.]
To move from $v$ to $w$, let us root the cube complex $X$ at $v$, and let $P$ be the corresponding PIP. Then $w$ corresponds to an order ideal $I$ of $P$; these are the hyperplanes we need to cross. 

\smallskip

\noindent
\emph{Cost ($\ell_1$) Metric}:
We simply cross the hyperplanes from $v$ to $w$  in non-decreasing order, with respect to the poset $I \subseteq P$: we first cross a minimal element $m_1 \in I$, then a minimal element $m_2 \in I-m_1$, and so on.

\smallskip

\noindent
\emph{Time ($\ell_\infty$) Metric}:
We first cross all minimal hyperplanes $M_1$ in $I$ simultaneously, then we cross all the minimal hyperplanes $M_2$ in $I-M_1$ simultaneously, and so on. This corresponds to Niblo and Reeves's \emph{normal cube path} \cite{NibloReeves}, where we cross the best available cube at each stage.
\end{proof}

These algorithms show how to move a CAT(0) robot optimally and automatically.

\section{Automation}

Driving in San Francisco, I get stuck behind a terrible driver. They are going extremely slowly, hesitating at every corner, stalling at every speed bump. When I finally lose patience and decide to pass them, they swerve wildly towards me; I  react quickly to avoid being hit. I turn to give the driver a nasty look, but I find there isn't one. 

What happens if you are injured by an automated, self-driving vehicle or robot designed by well-meaning scientists and technologists? When you live this close to Silicon Valley, the question is not just philosophical.

%

\section{Prototype: A robotic arm in a tunnel}\label{sec:prototype}

If we wish to apply Theorem \ref{thm:algorithm} to move an object optimally, our first hope is that the corresponding map of possibilities is a CAT(0) cubical complex. If this is true, we can prove it by choosing a convenient root and identifying the corresponding PIP. Tia Baker and Rika Yatchak pioneered this approach in their Master's theses \cite{ArdilaBakerYatchak}.

For concreteness, let us consider our robotic arm of length $n$ in a rectangular tunnel of height $1$. Baker and Yatchak found that the number of states of the configuration space is the term  $F_{n+1}$ of the Fibonacci sequence. This seemed like good news, until we realized that these numbers grow exponentially! The dimension of the map is $n/3$, and its combinatorial structure is enormous and intricate. We cannot navigate this map by brute force.

Fortunately, by running the bijection of Theorem \ref{thm:PIP} on enough examples, Baker and Yatchak discovered that this robot has a very nice PIP: a triangular wedge $T_n$ of a square grid with no inconsistent pairs, as shown in Figure \ref{fig:width1}. It is much simpler and only has about $n^2/4$ vertices. Indeed they proved that the map of possibilities of the robot $R_{1,n}$ is isomorphic to the cubical complex $X(T_n)$ corresponding to $T_n$. This implies that the map \textbf{is} CAT(0), \textbf{and} it allows us to use $T_n$ as a remote control to move the robot optimally.
\begin{center}
 \includegraphics[height=1.3in]{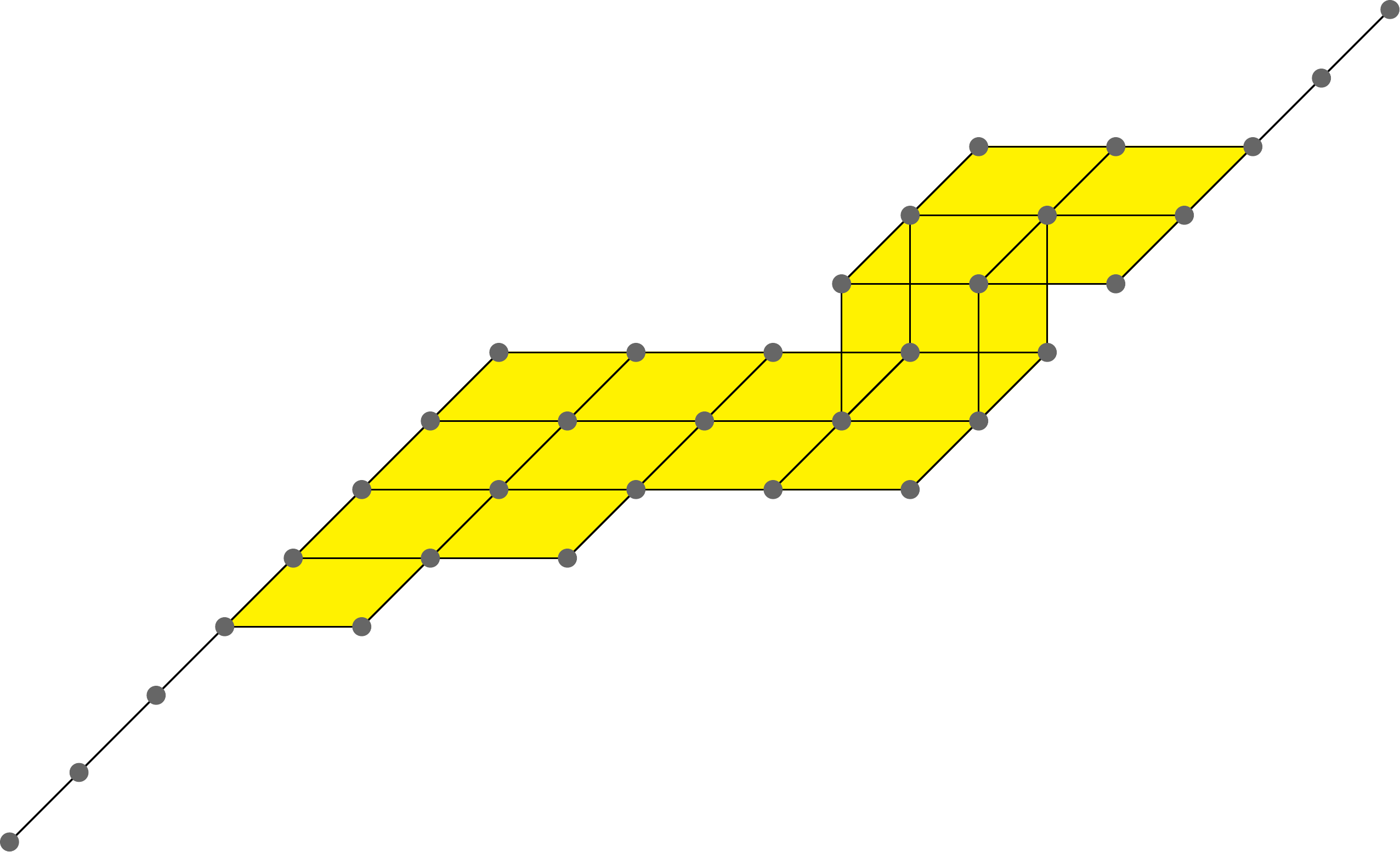}    \qquad
 \includegraphics[height=1.3in]{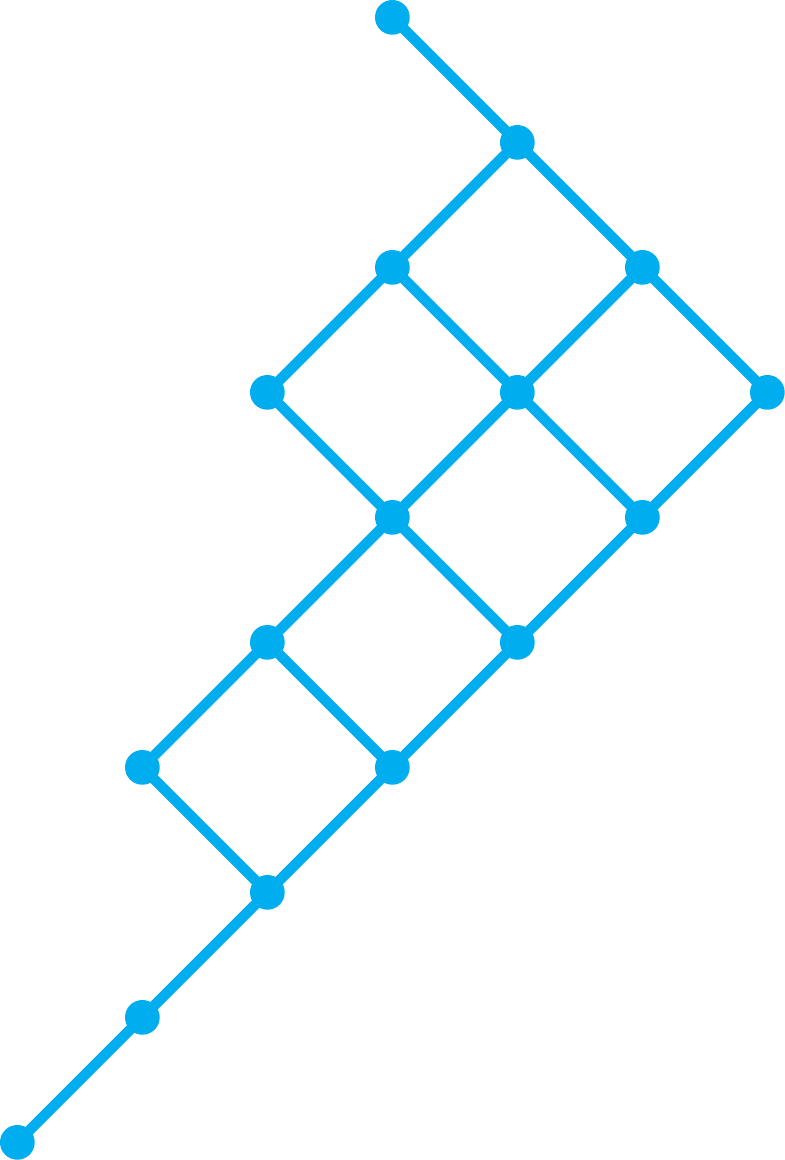}
\captionof{figure}{The map of the robot $R_{1,7}$ and its PIP $T_7$.\label{fig:width1}}
\end{center}

More generally, we have the following.

\begin{theorem}\label{th:CAT(0)2} \cite{ArdilaBakerYatchak, ArdilaBastidasCeballosGuo}
The configuration space of the robotic arm $R_{m,n}$ of length $n$ in a tunnel of height $m$ is a $\mathrm{CAT}(0)$ cubical complex. Therefore, we have an algorithm to move the arm optimally from any position to any other.
\end{theorem}

Naturally, as the height grows, the map becomes increasingly complex. 
After staring at many examples, getting stuck, and finally receiving a conclusive hint from the Pacific Ocean -- a piece of coral with a fractal-like structure -- we were able to describe the PIP of the robot $R_{m,n}$ for any $m$ and $n$. It is made of triangular flaps like the one in Figure \ref{fig:width1} recursively branching out in numerous directions.

This \emph{coral PIP} serves as a witness that the map of possibilities of the robotic arm $R_{m,n}$ is a CAT(0) cubical complex. It can also be programmed to serve as a remote control, to help the arm explore the tunnel.

\begin{center}
\includegraphics[height=2.1in]{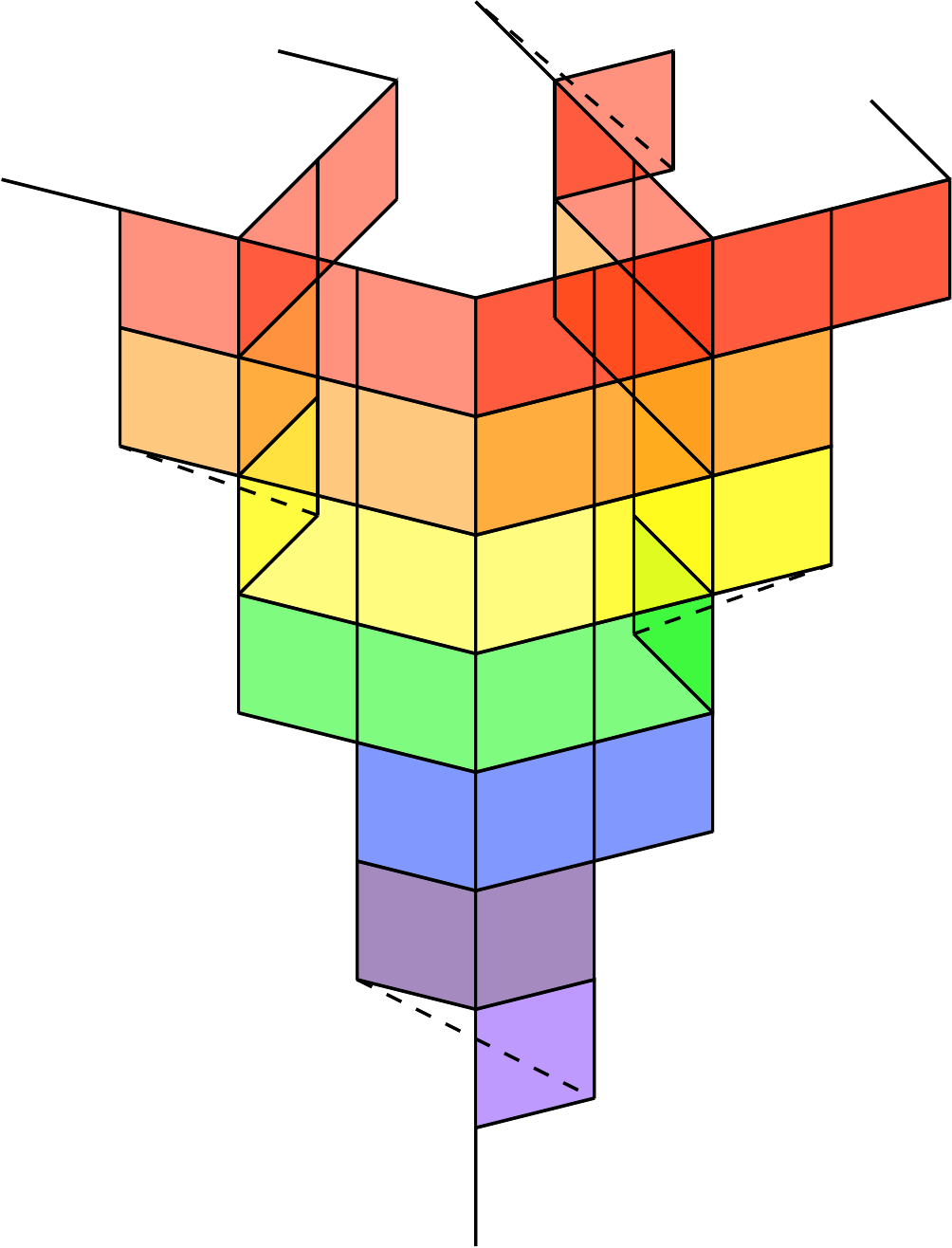} 
\captionof{figure}{The coral PIPs of the robot $R_{2,9}$, which contains the PIPs of $R_{2,1}, \ldots, R_{2,8}$, shown in different colors.}
\end{center}

\section{Implementation}

The algorithms to navigate a CAT(0) space optimally -- and hence move a CAT(0) robot, are described in~\cite{ArdilaBakerYatchak}. We have implemented them in Python for the robotic arm in a tunnel \cite{ArdilaBastidasCeballosGuo}. Given two states, the program outputs the distance between the two states in terms of cost ($\ell_1$) and time ($\ell_\infty$), and an animation moving the robot optimally between the two states. The downloadable code, instructions, and a sample animation are at \href{http://math.sfsu.edu/federico/robots.html}{\texttt{http://math.sfsu.edu/federico/robots.html}}.

With the goal to broaden access to these tools, I joined my collaborator C\'esar Ceballos, who led a week-long workshop for young robotics enthusiasts, as part of the Clubes de Ciencia de Colombia. This program invites Colombian researchers to design scientific activities for groups of students from public high schools and universities across the country.

We proposed some discrete models of robotic arms, and our students successfully built their maps of possibilities.
Extremely politely, they also pointed out that C\'esar and I really didn't know much about the mechanics of robots, and cleverly proposed several possible mechanisms. After the workshop, Arlys Asprilla implemented the design on CAD and built an initial prototype.

\begin{center}
  \includegraphics[height=2.1in]{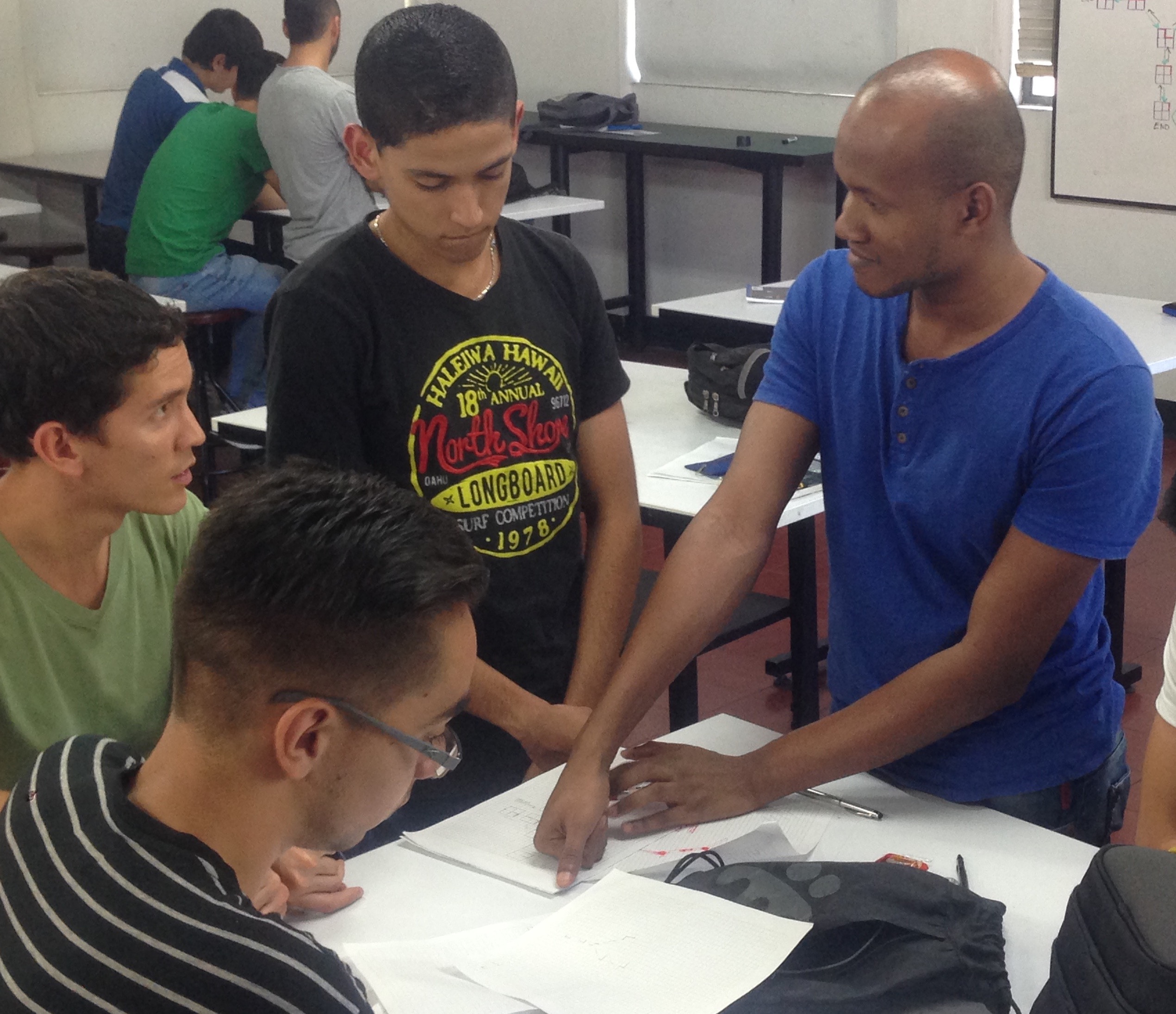}
  \includegraphics[height=2in]{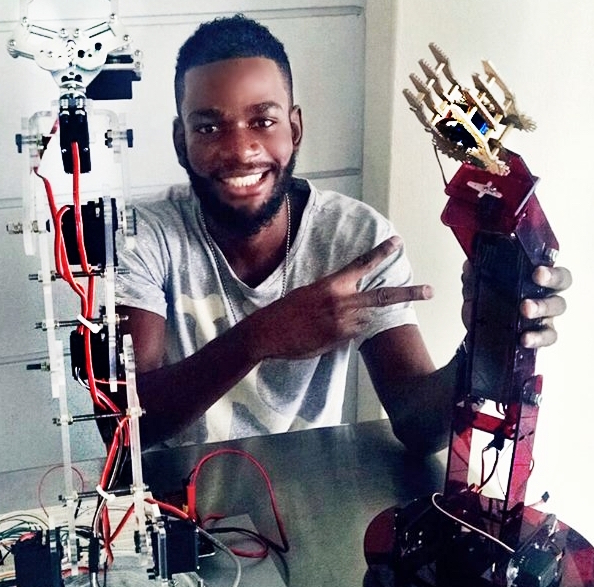} 
\captionof{figure}{\label{fig:clubesdeciencia}a. C\'esar Ceballos and students discuss configuration spaces during the Clubes de Ciencia de Colombia. b. Arlys Asprilla and one of his robotic arms.}
\end{center}

\section{Escuela de Rob\'otica del Choc\'o}

Arlys, his classmate Wolsey Rubio (on the right in Figure \ref{fig:clubesdeciencia}.a.), my partner May-Li Khoe, our friend Akil King, and I designed a similar workshop in Arlys and Wolsey's native Choc\'o. This region of the Colombian Pacific Coast  is one of the most biodiverse in the world, and also one of the most neglected historically by our government. We partnered with the Escuela de Rob\'otica del Choc\'o, led by Jimmy Garc\'{\i}a, which seeks to empower local youth to developg their scientific and technological skills, in order to address the problems faced by their communities.

\begin{center}
  \includegraphics[height=2in]{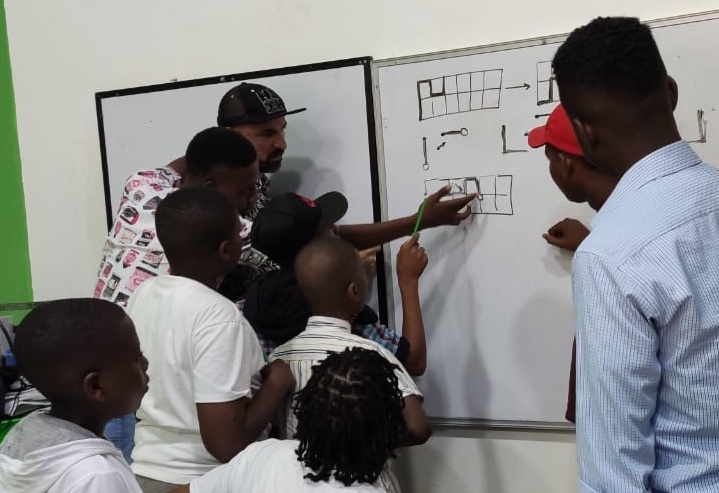}
  \includegraphics[height=2.4in]{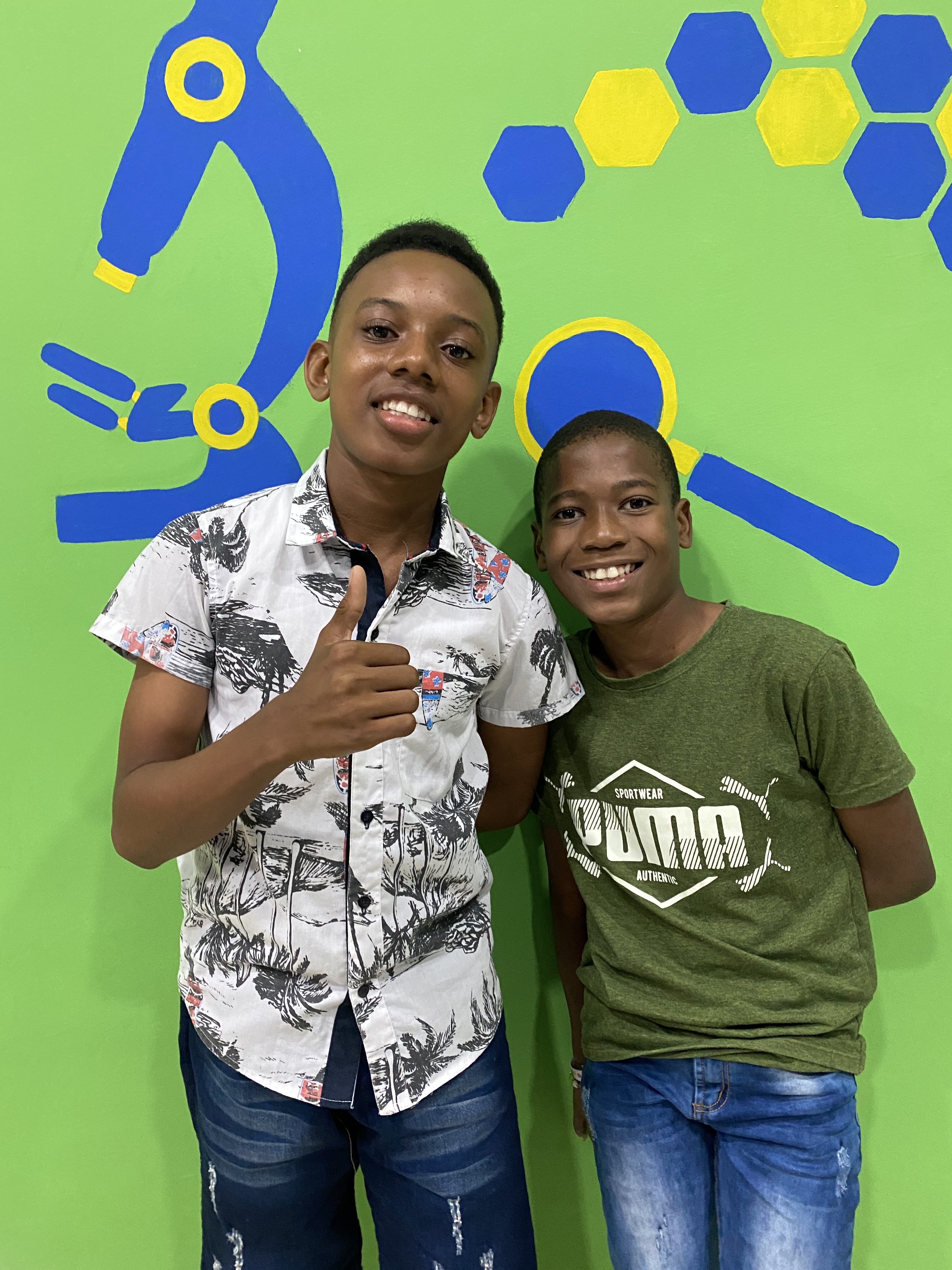} 
\captionof{figure}{a. At the Escuela de Rob\'otica del Choc\'o. b. Deison Rivas and Juan David Cuenta.}
\end{center}

At the end of the workshop, we asked the students: What robot do you \textbf{really} want to design?
Deison Rivas wants to build a firefighter robot; it will quickly and safely go in and out of houses --  traditionally made of wood -- and put out the fires that have razed entire city blocks in Quibd\'o in the past. 
Juan David Cuenta wants to design an agile rescue robot; it will help people stuck under the frequent landslides caused by illegal mining operations and by heavy rainfalls on the roads.

This theoretical exercise on robotic optimization immediately took on new meaning, thanks to the wisdom of these young people.

\section{Our role as educators}

As I put the finishing touches on this paper, I remember a book I inherited from my mom when she passed away ten years ago. Browsing through it I find this passage, highlighted by her:

\begin{quote}
The saddest achievement of our educational system is to produce [...]
the most disastrous person of our society: the creative scientist who is at the same time enslaved to the military or industrial apparatus; someone who makes contributions, but has no interest in the way they will be used.

-- Estanislao Zuleta \cite{Zuleta}
\end{quote}

\section{What does it mean to do math ethically?}

Six years ago, my student Brian Cruz asked me whether mathematicians have an ethical code, similar to the Hippocratic Oaths adopted by physicians. More than two decades into my mathematical career, I had never thought or heard of this specific suggestion.

Thanks to Brian, I did some research, gathered some resources with the help of my students\footnote{These resources are available at \href{http://math.sfsu.edu/federico/ethicsinmath.html}{\texttt{http://math.sfsu.edu/federico/ethicsinmath.html}}}, and I now devote one day of each semester to discuss this question with them. 
Posing the question to them is surely more important than proposing an answer:

\smallskip

\noindent
\emph{\textbf{Writing assignment.} What does ``doing mathematics ethically" mean to you? This question is an invitation to recognize the power you carry as a mathematician, and the privilege and responsibility that comes with it. When you enter a scientific career, you do not leave yourself at the door. You can choose how to use that power. My hope is that you will always continue to think about this in your work.}

\section{Acknowledgments} 

Este art\'{\i}culo est\'a dedicado a mi mam\'a, Amparo Mantilla. Una de sus incontables ense\~nanzas es la importancia de pensar cr\'{\i}ticamente en el impacto de la ciencia en la sociedad.

I would like to extend my sincere gratitude to the students and colleagues who have collaborated with me on this research: Arlys Asprilla, C\'esar Ceballos, Hanner Bastidas, John Guo, Matthew Bland, Maxime Pouokam, Megan Owen, Rika Yatchak, Seth Sullivant, and Tia Baker.

\noindent
\begin{center}
\includegraphics[height=0.8in]{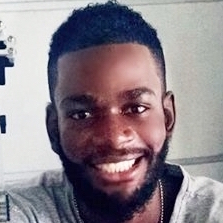}
\includegraphics[height=0.8in]{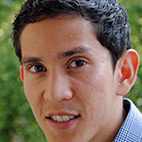}
\includegraphics[height=0.8in]{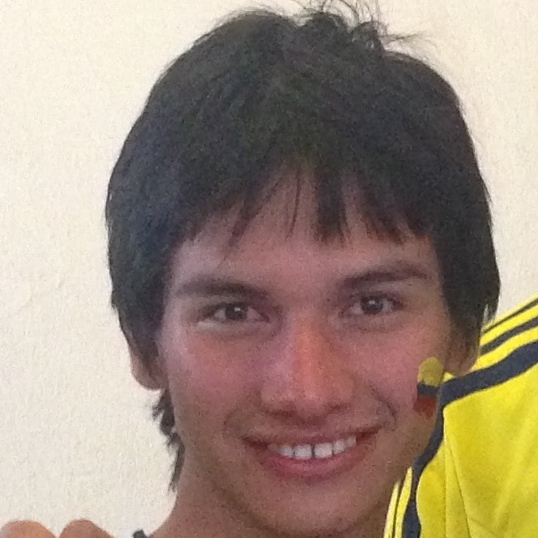}
\includegraphics[height=0.8in]{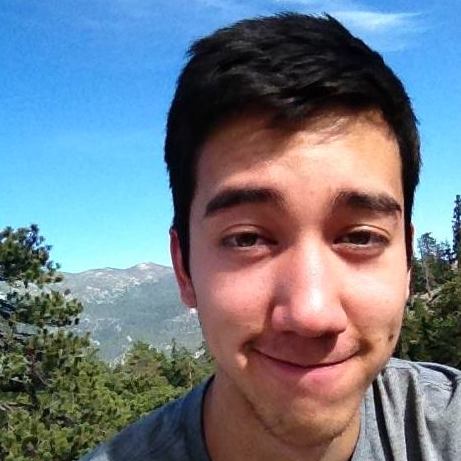}
\includegraphics[height=0.8in]{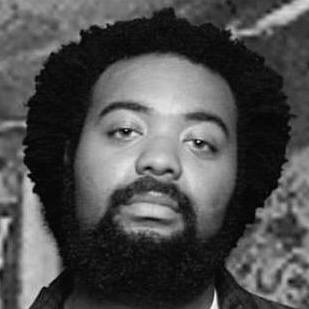}
\includegraphics[height=0.8in]{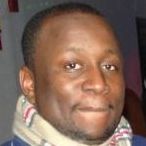}
\includegraphics[height=0.8in]{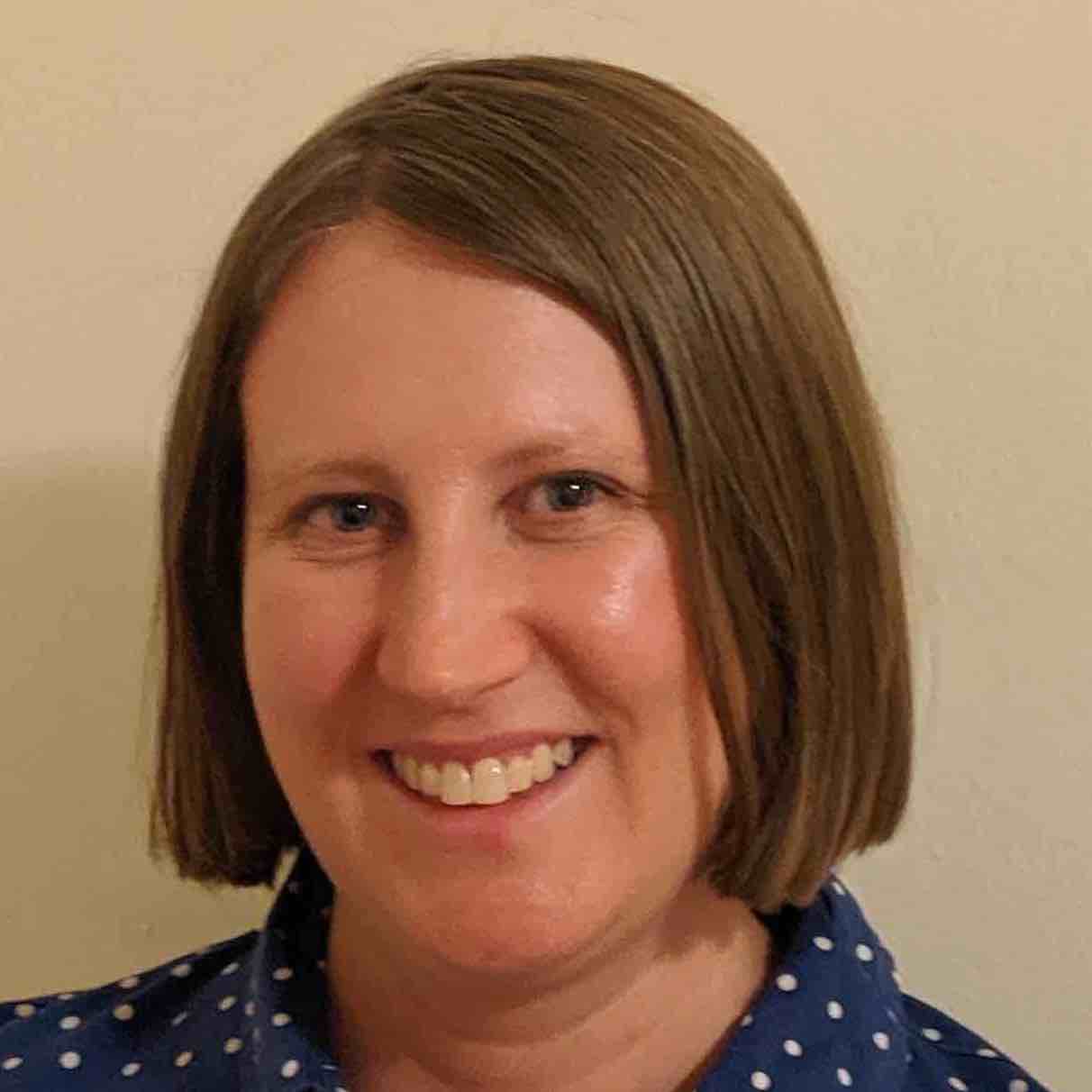} 
\includegraphics[height=0.8in]{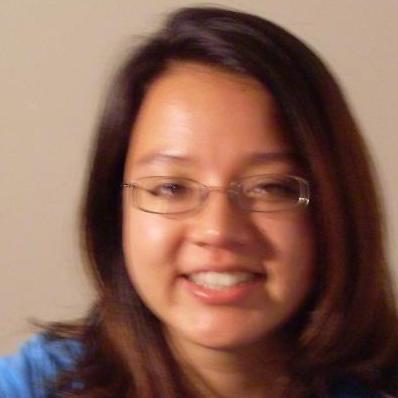}
\includegraphics[height=0.8in]{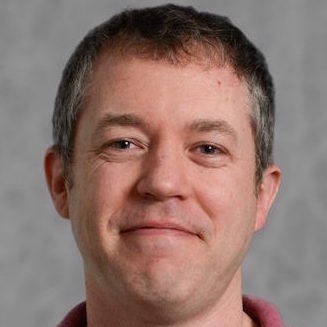} 
\includegraphics[height=0.8in]{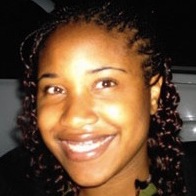} 
\end{center}

\noindent
I am also very grateful to the many people who encouraged me and helped me figure out how to tell this story.

\footnotesize

\makeatletter
\renewcommand{\@biblabel}[1]{\hfil#1.}
\makeatother
\bibliographystyle{myamsplain}
\bibliography{biblio}

\end{multicols}

\vfill

%

\vfill

\end{document}